\documentclass[a4paper,12 pt]{article}

\usepackage{amsmath}
\usepackage{amssymb}
\usepackage[english]{babel}
\usepackage{}
\usepackage[ansinew]{inputenc}
\usepackage{graphicx}
\usepackage{color}
\usepackage[colorlinks]{hyperref}
\newtheorem{thm}{Theorem}[section]

\newtheorem{lemm}[thm]{Lemma}
\newtheorem{remark}[thm]{Remark}

\newtheorem{notation}[thm]{Notation}

\newtheorem{defin}[thm]{Definition}
\newtheorem{defs}[thm]{Definitions}
 \numberwithin{equation}{section}
\def\R{\mathbb R}

\def\epsilon{\varepsilon}
\def\o{\overline\Omega}
\def\ro{\R\times\o}

\def\te{\tilde{e}}
\def\ds{\displaystyle}
\allowdisplaybreaks
\begin{document}
\title{ \textbf{Min-Max formul{\ae} for the speeds of pulsating travelling fronts in periodic excitable media} }
\author{Mohammad El Smaily\footnote{E-mail: elsmaily@math.ubc.ca} \\
\begin{small} Universit\'e Aix-Marseille III, LATP, Facult\'e des
Sciences et
Techniques,\end{small}\\
 \begin{small}Avenue Escadrille Normandie-Niemen, F-13397 Marseille
Cedex 20, France.\end{small}}
%\date{27 March 2008}
\maketitle

 \paragraph{Abstract.} This paper is concerned with some nonlinear propagation phenomena for reaction-advection-diffusion equations in a periodic framework. It deals with travelling wave solutions of the equation $$u_t =\nabla\cdot(A(z)\nabla u)\;+q(z)\cdot\nabla u+\,f(z,u),~~ t\,\in\,\mathbb{R},\;z\,\in\,\Omega,$$
 propagating with a speed $c.$ In the case of a ``combustion'' nonlinearity, the speed $c$ exists and it is unique, while the front $u$ is unique up to a translation in $t.$ We give a $\min-\max$ and a $\max-\min$ formula for this speed $c.$ On the other hand, in the case of a ``ZFK'' or a ``KPP'' nonlinearity, there exists a minimal speed of propagation $c^{*}.$ In this situation, we give a $\min-\max$ formula for $c^{*}.$ Finally, we apply this $\min-\max$ formula to prove a variational formula involving eigenvalue problems for the minimal speed $c^{*}$ in the ``KPP'' case.

\tableofcontents
\section{Introduction and main results}
\subsection{A description of the periodic framework}
The goal of this paper is to give some variational formul{\ae} for
the speeds of pulsating travelling fronts corresponding to
reaction-diffusion-advection equations set in a heterogenous
periodic framework. In fact, many works, such as Hamel \cite{Hamel},
Heinze, Papanicolaou, Stevens \cite{HPS}, and Volpert, Volpert, Volpert \cite{3volpert} treated this problem
in simplified situations and under more strict assumptions. In this paper, we treat the problem in the most general periodic framework. We are concerned with equations of the type

 \begin{eqnarray}\label{heteq}
    \left\{
      \begin{array}{l}
 u_t =\nabla\cdot(A(z)\nabla u)\;+q(z)\cdot\nabla u+\,f(z,u),\quad t\,\in\,\mathbb{R},\;z\,\in\,\Omega,
 \vspace{4 pt}\\
        \nu \cdot A\;\nabla u(t,z) =0,\quad
        t\,\in\,\mathbb{R},\; z\,\in\,\partial\Omega,
      \end{array}
    \right.
\end{eqnarray}
where $\nu(z)$ is the unit outward normal on $\partial\Omega$ at the
point $z.$ In this context, let us detail the mathematical
description of the heterogeneous setting.

 Concerning the domain, let $N\geq\,1$ be the
space dimension, and let $d$ be an integer so that
$1\leq\,d\leq\,N.$ For an element
$z=(x_1,x_2,\cdots,x_d,x_{d+1},\cdots,x_N)\in\,\mathbb{R}^{N},$ we
denote by $x\,=\,(x_1,x_2,\cdots,x_d)$ and by $\;
y=(x_{d+1},\cdots,x_N)$ the two tuples so that $z=(x,y).$ Let
$L_1,\cdots,L_d$ be $d$ positive real numbers, and let $\;\Omega\;$
be a $C^3$ non empty connected open subset of $\mathbb{R}^N$
satisfying
\begin{eqnarray}\label{comega}
    \left\{
      \begin{array}{l}
        \exists\,R\geq0\,;\forall\,(x,y)\,\in\,\Omega,\,|y|\,\leq\,R, \vspace{3pt} \\
        \forall\,(k_1,\cdots,k_d)\;\in\;L_1\mathbb{Z}\times\cdots\,\times L_d\mathbb{Z},
        \quad\displaystyle{\Omega\;=\;\Omega+\sum^{d}_{k=1}k_ie_i},
      \end{array}
    \right.
\end{eqnarray}
where $\;(e_i)_{1\leq i\leq N}\;$ is the canonical basis of
$\mathbb{R}^N.$

 As $d\geq1,$ one notes that the set $\Omega$ satisfying (\ref{comega}) is unbounded. We have many archetypes of such a domain. The case of the whole space $\mathbb{R}^{N}$ corresponds to $d =
N,$ where $L_1,\ldots,L_N$ are any positive numbers. The case of the
whole space $\mathbb{R}^{N}$ with a periodic array of holes can also
be considered. The case $d = 1$ corresponds to domains which have only
one unbounded dimension, namely infinite cylinders which may be
straight or have oscillating periodic boundaries, and which may or
may not have periodic perforations. The case $2\leq d \leq N-1$ corresponds
to infinite slabs.

In this periodic situation, we give the following definitions:
\begin{defin}[Periodicity cell] The set $$C\;=\;\{\,(x,y)\in\Omega
;\; x_{1}\in(0,L_1),\cdots,x_{d}\in(0,L_d)\}$$ is called the
periodicity cell of $\Omega.$
\end{defin}
\begin{defin}[$L$-periodic fields]  A field
$w:\Omega\rightarrow\,\mathbb{R}^N$ is said to be $L$-periodic with
respect to $x$
 if $$w(x_1+k_1,\cdots,x_d+k_d\,,y)\;=\;w(x_1,\cdots,x_d,y)$$ almost everywhere in
$ \Omega,$ and for all
$\displaystyle{k=(k_1,\cdots\,,k_d)\in\prod^{d}_{i=1}L_i\mathbb{Z}}.$
\end{defin}

Let us now detail the assumptions concerning the coefficients in
 (\ref{heteq}). First, the diffusion matrix $A(x,y)=(A_{ij}(x,y))_{1\leq i,j\leq
 N}$ is a symmetric $C^{2,\delta}(\,\overline{\Omega}\,)$ (with $\delta
\,>\,0$) matrix field satisfying
\begin{eqnarray}\label{cA}
    \left\{
      \begin{array}{l}
        A\; \hbox{is $L-$periodic with respect to}\;x,  \vspace{4 pt}\\
        \exists\,0<\alpha_1\leq\alpha_2;\forall(x,y)\;\in\;\Omega,\forall\,\xi\,\in\,\mathbb{R}^N, \vspace{4 pt}\\
        \displaystyle{\alpha_1|\xi|^2 \;\leq\;\sum_{1\leq i,j\leq N} \,A_{ij}(x,y)\xi_i\xi_j\,\;\leq\alpha_2|\xi|^2.}
      \end{array}
    \right.
\end{eqnarray}

The underlying advection $q(x,y)=(q_1(x,y),\cdots,q_N(x,y))$  is a
$C^{1,\delta}(\overline{\Omega})$ (with $\delta>0$) vector field
satisfying:
\begin{equation}\label{cq}
    \left\{
      \begin{array}{ll}
        q\;\hbox{is $L-$periodic with respect to }\;x, & \hbox{} \vspace{4 pt} \\
        \nabla\cdot q\,=0\quad\hbox{in}\; \overline{\,\Omega\,}, \vspace{4 pt} \\
         q\cdot\nu\,=0\quad \hbox{on}\;\partial\Omega\,, \vspace{4 pt} \\
        \forall\,1\leq\,i\,\leq\,d, \quad\displaystyle{\int_{C}q_i\;dx\,dy \,=\,0}\hbox{.}
      \end{array}
    \right.
\end{equation}

Lastly, let $f=f(x,y,u)$ be a function defined in
$\overline{\Omega}\times\mathbb{R}$ such that
\begin{eqnarray}\label{f1}
    \left\{
      \begin{array}{ll}
f \;\hbox{is globally Lipschitz-continuous in}\,
\;\overline{\Omega}\times\mathbb{R}, \vspace{4 pt}\\
\displaystyle{\forall\,(x,y)\in\,\overline{\Omega},\,\forall\,s\in(-\infty,0]\cup[1,+\infty),\,f(s,x,y)=0,} \vspace{4 pt}\\
 \exists
\,\rho\in\,(0,1),\;\forall(x,y)\,\in\overline{\Omega},\;\displaystyle{\forall\,
1-\rho\leq\,s\,\leq s'\,\leq\,1,}\;\vspace{3pt}\\
\displaystyle{f(x,y,s)\;\geq\,f(x,y,s')}.
  \end{array}
    \right.
\end{eqnarray}

One assumes that
\begin{equation}\label{f2}
f\;\hbox{is $L-$periodic with respect to}\; x.
\end{equation}

Moreover, the function $f$ is assumed to be of one of the following two
types: either
\begin{eqnarray}\label{combustion}
\left\{
\begin{array}{ll}
\exists\; \theta\,\in\,(0,1),\;
\forall(x,y)\in\overline{\Omega},\;\forall\, s\in[0,\theta],\;
f(x,y,s)=0,\vspace{3pt}\\
\forall\,s\in\,(\theta,1),\; \exists
\,(x,y)\in\overline{\Omega}\;\;\hbox{such that}\;f(x,y,s)\,>\,0
\hbox{,}
\end{array}
\right.
\end{eqnarray}

or
\begin{eqnarray}\label{ZFK}
    \left\{
      \begin{array}{ll}
\exists\;\delta>0,\;\hbox{the restriction of $f$ to $\overline{\Omega}\,\times\,[0,1]$ is of class}\;C^{1,\,\delta},  \vspace{3pt}\\
 \forall\,s\in\,(0,1),\; \exists \,(x,y)\in\overline{\Omega}\;\;\hbox{such that}\;f(x,y,s)\,>\,0  \hbox{.} \\
\end{array}
    \right.
\end{eqnarray}
\begin{defs}
A nonlinearity $f$ satisfying (\ref{f1}), (\ref{f2}) and
(\ref{combustion}) is called a ``combustion'' nonlinearity. The
value $\theta$ is called the ignition temperature.\\
A nonlinearity $f$ satisfying (\ref{f1}), (\ref{f2}), and
(\ref{ZFK}) is called a ``ZFK'' (for Zeldovich-Frank- Kamenetskii) nonlinearity.\\
If $f$ is a ``ZFK'' nonlinearity that satisfies
\begin{equation}\label{nondegeneracy}
\displaystyle{f'_{u}(x,y,0)=\lim_{u\rightarrow0^+}\;f(x,y,u)/u\;>0,}
\end{equation}
with the additional assumption
\begin{equation}\label{KPP condition}
\forall\, (x,y,s)\in\overline{\Omega}\times(0,1),\quad
0\,<\,f(x,y,s)\,\leq\,f'_u(x,y,0)\times\,s ,
\end{equation}
then $f$ is called a ``KPP''(for Kolmogorov, Petrovsky, and Piskunov, see \cite{KPP}) nonlinearity.
\end{defs}

The simplest examples of ``combustion'' and ``ZFK''
nonlinearities are when $f(x,y,u)=f(u)$ where: either
\begin{eqnarray}\label{combustion hom}
\left\{
  \begin{array}{ll}
f\;\hbox{is Lipschitz-continuous in}\; \R,\;\vspace{3pt}\\
\exists\,\theta\in\,(0,1),\;f(s)=0\;\hbox{for
all}\;s\in\,(-\infty,\theta]\cup[1,+\infty),\vspace{3pt}\\
f(s)>0\;\hbox{for all}
\;s\in(\theta,1),\vspace{3pt}\\
\exists \;\rho\in(0,1-\theta),\quad f\;\hbox{is non-increasing on}
\;[1-\rho,1],
  \end{array}
\right.
\end{eqnarray}

or

\begin{eqnarray}\label{ZFK homogenous}
\left\{
\begin{array}{ll}
f\; \hbox{is defined on $\mathbb{R},$}\;f\equiv0\;\hbox{in}\;\mathbb{R}\setminus{(0,1)},\vspace{3pt}\\
 \exists\;\delta>0,\; \hbox{the restriction of $f$ on the
interval $[0,1]$ is } C^{1,\delta}([0,1]),\vspace{3pt}\\
f(0)=f(1)=0,\;\hbox{and}\;f(s)\,>0\;\hbox{ for all}\;s\in(0,1),\vspace{3pt}\\
\exists \;\rho>0,\quad f\;\hbox{is non-increasing on} \;[1-\rho,1].
\end{array}
 \right.
\end{eqnarray}

If $f=f(u)$ satisfies (\ref{combustion hom}), then it is a
homogeneous ``combustion'' nonlinearity. On the other hand, a
nonlinearity $f=f(u)$ that satisfies (\ref{ZFK homogenous}) is
homogeneous of the ``ZFK'' type. Moreover, a KPP homogeneous
nonlinearity is a function $f=f(u)$ that satisfies (\ref{ZFK
homogenous}) with the additional assumption
\begin{equation}\label{kPP condition for f=f(u)}
\forall\,s\in(0,1),\;0<f(s)\leq f'(0)\,s.
\end{equation}

As typical examples of nonlinear heterogeneous sources satisfying
(\ref{f1}-\ref{f2}) and either (\ref{combustion}) or (\ref{ZFK}),
one can consider the functions of the type $$f(x,y,u)=h(x,y)\,g(u),$$ where
$h$ is a globally Lipschitz-continuous, positive, bounded, and
$L-$periodic with respect to $x$ function defined in
$\overline{\Omega},$ and $g$ is a function satisfying either
(\ref{combustion hom}) or (\ref{ZFK homogenous}).

\begin{defin}[Pulsating fronts and speed of propagation]

Let $e=(e^1,\cdots,e^d)$ be an arbitrarily given unit vector in
$\mathbb{R}^d.$ A function $u\,=\,u(t,x,y)$ is called a pulsating
travelling front propagating in the direction of $-e$ with an
effective speed $c\,\neq\,0,$ if $u$ is a classical solution of:
\begin{eqnarray}\label{front}
    \left\{
      \begin{array}{ll}
u_t =\nabla\cdot(A(x,y)\nabla u)+q(x,y)\cdot\nabla u+\,f(x,y,u),\;
t\,\in\,\mathbb{R},\;(x,y)\,\in\,\Omega,
 \vspace{4 pt}\\
        \nu \cdot A\;\nabla u(t,x,y) =0,\;
        t\,\in\,\mathbb{R},\;(x,y)\,\in\,\partial\Omega, \vspace{4 pt}\\
 \displaystyle{\forall\, k\in\prod^{d}_{i=1}L_i\mathbb{Z},\, \forall\,(t,x,y)\,\in\,\mathbb{R}\times\overline{\Omega}},
\;\displaystyle{u(t+\frac{k\cdot e}{c},x,y)=u(t,x+k,y)}  \hbox{,} \vspace{4 pt} \\
         \displaystyle{\lim_{x\cdot e\rightarrow-\infty}u(t,x,y)=0,\;\hbox{and}\; \lim_{x\cdot e\rightarrow +\infty} u(t,x,y)=1}  \hbox{,}
          \vspace{4 pt}\\
          0\leq u\leq1,
      \end{array}
    \right.
\end{eqnarray}
 where the above limits hold locally in $t$ and uniformly in $y$ and in the
directions of $\mathbb{R}^d$ which are orthogonal to $e$ .
\end{defin}

 Several works were concerned with pulsating travelling fronts in periodic media (see \cite{BH1}, \cite{BHN1}, \cite{KKTS}, \cite{PapXin}, \cite{shigesada Kawsaki 1}, \cite{shigesada Kawsaki 2}, and \cite{JXin}).

In the general periodic framework, we recall two essential known
results and then we move to our main results.
\begin{thm}[Berestycki, Hamel \cite{BH1}]\label{Berest Hamel combustion Thm}
Let $\Omega$ be a domain satisfying (\ref{comega}), let $e$ be any
unit vector of $\mathbb{R}^{d}$ and let $f$ be a nonlinearity
satisfying (\ref{f1}-\ref{f2}) and (\ref{combustion}). Assume,
furthermore, that $A$ and $q$ satisfy (\ref{cA}) and (\ref{cq})
respectively. Then, there exists a classical solution $(c,u)$ of
(\ref{front}). Moreover, the speed $c$ is positive and unique while
the function $u=u(t,x,y)$ is increasing in $t$ and it is unique up
to a translation. Precisely, if $(c^{1},u^{1})$ and $(c^{2}, u^{2})$
are two classical solutions of (\ref{front}), then $c^{1}=c^{2}$ and
there exists $h\in\mathbb{R}$ such that
$u^{1}(t,x,y)=u^{2}(t+h,x,y)$ for all
$(t,x,y)\in\mathbb{R}\times\overline{\Omega}.$
\end{thm}

In a periodic framework, Theorem \ref{Berest Hamel combustion Thm}
yields the existence of a pulsating travelling front in the case of a
``combustion'' nonlinearity with an ignition temperature $\theta.$
It implies, also, the uniqueness of the speed and of the profile of
$u.$ For  ``ZFK'' nonlinearities, we have
\begin{thm}[Berestycki, Hamel \cite{BH1}]\label{Berest Hamel ZFK Thm}
Let $\Omega$ be a domain satisfying (\ref{comega}), let $e$ be any
unit vector in $\mathbb{R}^{d}$ and let $f$ be a nonlinearity
satisfying (\ref{f1}-\ref{f2}) and (\ref{ZFK}). Assume, furthermore,
that $A$ and $q$ satisfy (\ref{cA}) and (\ref{cq}) respectively.
Then, there exists $\displaystyle{c^{*}_{\Omega,A,q,f}(e)>0}$ such that
the problem (\ref{front}) has no solution $(c,u)$ such that $u_t>0$ in $\ro$ if
$c<\displaystyle{c^{*}_{\Omega,A,q,f}(e)}$ while, for each
$c\geq\displaystyle{c^{*}_{\Omega,A,q,f}(e)},$ it has a solution
$(c,u)$ such that $u$ is increasing in $t.$
\end{thm}

In fact, the existence and the monotonicity of a solution $u^{*}=u^{*}(t,x,y)$ of (\ref{front}) for $c=\displaystyle{c^{*}_{\Omega,A,q,f}(e)>0}$ holds by approaching the ``ZFK'' nonlinearity $f$ by a sequence of combustion nonlinearities $\left(f_\theta\right)_\theta$ such that $f_\theta\rightarrow f$ uniformly in $\ro$ as $\theta \searrow 0^{+}$ (see more details in step 2 of the proof of formula (\ref{min max ZFK }) below, section \ref{section pf of min max zfk}). It follows, from Theorem \ref{Berest Hamel combustion Thm}, that for each $\theta>0,$ there exists a solution $(c_\theta,u_\theta)$ of (\ref{front}) with the nonlinearity $f_\theta$ such that $u_\theta$ is increasing with respect to $t.$ From parabolic estimates, the functions $u_\theta,$ converge up to a subsequence, to a function $u^{*}$ in $C^{2}_{loc}(\ro)$ as $\theta\rightarrow0^{+}.$ Moreover, Lemmas 6.1 and 6.2 in \cite{BH1} yield the existence of a constant $c^{*}(e)=\displaystyle{c^{*}_{\Omega,A,q,f}(e)}>0$ such that $c_\theta\nearrow c^{*}(e)$ as $\theta\searrow0.$ Hence, the couple $(c^{*}(e),u^{*})$ becomes a classical solution of (\ref{front}) with the nonlinearity $f$ and one gets
 that $u^{*}$ is nondecreasing with respect to $t$ as a limit of the increasing functions $u_\theta.$ Finally, one applies the strong parabolic maximum principle and Hopf lemma to get that $w$ is positive in $\ro.$ In other words, $u^{*}$ is increasing in $t\in\R.$ Actually, in the ``ZFK'' case, under the additional non-degeneracy assumption (\ref{nondegeneracy}), it is known that any pulsating front with speed $c$ is increasing in time and
 $c\geq c^{*}(e)$ (see \cite{BH1}).

The value $\displaystyle{c^{*}_{\Omega,A,q,f}(e)}$ which appears in
Theorem \ref{Berest Hamel ZFK Thm} is called the minimal speed of
propagation of the pulsating travelling fronts propagating in the
direction $-e$ (satisfying the reaction-advection-diffusion problem
(\ref{front})).

We mention that the uniqueness of the pulsating
travelling fronts, up to shifts in time, for each
$c\geq\displaystyle{c^{*}_{\Omega,A,q,f}(e)},$ has been proved recently by Hamel and Roques \cite{Hamel Roques} for ``KPP'' nonlinearities. On the other hand, a variational formula for
the minimal speed of propagation
$\displaystyle{c^{*}_{\Omega,A,q,f}(e)},$ in the case of a KPP
nonlinearity, was proved in Berestycki, Hamel, Nadirashvili
\cite{BHN1}. This formula involves eigenvalue problems and gives the
value of the minimal speed in terms of the domain
$\Omega$ and in terms of the coefficients $A,\;q,$ and $f$ appearing
in problem (\ref{front}). The asymptotic behaviors and the
variations of the minimal speed of propagation, as a function of the
diffusion, advection and reaction factors and as a function of the
periodicity parameters, were widely studied in Berestycki, Hamel, Nadirashvili \cite{BHN2}, El~Smaily \cite{El
Smaily}, El Smaily, Hamel, Roques \cite{EHR}, Heinze \cite{Heinze advection}, Ryzhik, Zlato\v{s} \cite{Zlatos Ryzhik}, and Zlato\v{s} \cite{zlatos speed up}.
\subsection{Main results}
In the periodic framework, having (in Theorems \ref{Berest Hamel
combustion Thm} and \ref{Berest Hamel ZFK Thm}) the existence
results and some qualitative properties of the pulsating travelling
fronts propagating in the direction of a fixed unit vector
$-e\in\mathbb{R}^{d},$ we search a variational formula for the unique
speed of propagation $c=c(e)$ whenever $f$ is of the ``combustion''
type, and for the minimal speed
$c^{*}=\displaystyle{c^{*}_{\Omega,A,q,f}(e)}$ whenever $f$ is of
the ``ZFK'' or the ``KPP'' type. We will answer the above
investigations in the following theorem, but before this, we
introduce the following
\begin{notation}\label{F(phi)}
For each function $\phi=\phi(s,x,y)$ in
$C^{1,\delta}(\mathbb{R}\times\overline{\Omega})$ {\rm(}for some
$\delta\in[0,1)${\rm )}, let
$$F[\phi]:=\nabla_{x,y}\cdot(A\nabla_{x,y}\phi)+(\tilde{e}A\tilde{e})\phi_{ss}+\nabla_{x,y}\cdot(A\tilde{e}\phi_s)+\partial_s(\tilde{e}A\nabla_{x,y}\phi)~\hbox{in}~\mathcal{D'}(\mathbb{R}\times\Omega),$$
where $\tilde{e}=(e,0,\cdots,0)\in\mathbb{R}^{N}$ and $e$ denotes a unit vector of $\R^d.$
\end{notation}
The first main result deals with the ``combustion'' case.
\begin{thm}\label{thm min max max min comb}
 Let $e$ a unit vector of $\mathbb{R}^{d}.$ Assume that $\Omega$ is
a domain satisfying (\ref{comega}) and $f$ is a nonlinear source satisfying (\ref{f1}) and (\ref{f2}). Assume furthermore that
$A$ and $q$ satisfy (\ref{cA}) and (\ref{cq}) respectively. Consider
the set of functions
\begin{eqnarray*}
\begin{array}{lll}
E&=&\left\{\varphi=\varphi(s,x,y),~\varphi~\hbox{is of class}~C^{1,\mu}(\mathbb{R}\times\overline{\Omega})~\hbox{for
each}~\mu\in[0,1),\right.\vspace{4pt}\\
&&
~F[\varphi]\in
C(\mathbb{R}\times\overline{\Omega}),\,\varphi~\hbox{is}\hbox{ $L-$periodic with respect to
$x,$}~\varphi_{s}>0~\vspace{4pt}\\
&&\hbox{in}~\mathbb{R}\times\overline{\Omega},~\varphi(-\infty,.,.)=0,~\varphi(+\infty,.,.)=1~ \hbox{uniformly
in}~\overline{\Omega},~\hbox{and}~~\vspace{4pt}\\
&&\left. \ds{\nu\cdot
A(\nabla_{x,y}\varphi+\tilde{e}\varphi_{s})=0\hbox{ on }\mathbb{R}\times\partial{\Omega}}\right\}.
\end{array}
\end{eqnarray*}
For each $\varphi\in E,$ we define the function $R\varphi\in
C(\mathbb{R}\times\overline{\Omega})$ as, for all
$(s,x,y)\in\mathbb{R}\times\overline{\Omega},$

\begin{eqnarray*}
\begin{array}{rl}
R\,\varphi(s,x,y)=&
\displaystyle{\frac{F[\varphi](s,x,y)+\,q\cdot\nabla_{x,y}\varphi(s,x,y)+f(x,y,\varphi)}{\partial_{s}\varphi(s,x,y)}+q(x,y)\cdot\te}.
\end{array}
\end{eqnarray*}
 If $f$ is a nonlinearity of ``combustion'' type satisfying
 (\ref{combustion}), then the unique speed $c(e)$
that corresponds to problem (\ref{front}) is given by
\begin{eqnarray}\label{min max combustion}
c(e)= \displaystyle{\min_{\varphi\in
E}\sup_{(s,x,y)\in\mathbb{R}\times\overline{\Omega}}R\,\varphi(s,x,y)},
\end{eqnarray}
\begin{eqnarray}\label{max min combustion}
c(e)&= \displaystyle{\max_{\varphi\in
E}\inf_{(s,x,y)\in\mathbb{R}\times\overline{\Omega}}R\,\varphi(s,x,y).}
\end{eqnarray}
Furthermore, the $\min$ in (\ref{min max combustion}) and the $\max$ in
(\ref{max min combustion}) are attained by, and only by, the function
$\phi(s,x,y)=u\left(\frac{s-x\cdot e}{c(e)},x,y\right)$ and its shifts $\phi(s+\tau,x,y)$ for any $\tau\in\R,$ where $u$ is
the solution of (\ref{front}) with a speed $c(e)$ (whose
existence and uniqueness up to a translation in $t$ follow from
Theorem \ref{Berest Hamel combustion Thm}).
\end{thm}
The second result is concerned with ``ZFK'' nonlinearities.
\begin{thm} \label{thm min max ZFK}
Under the same notations of Theorem \ref{thm min max max min comb},
if $f$ is a nonlinearity of ``ZFK'' type satisfying (\ref{f1}-\ref{f2}) and (\ref{ZFK}), then
 the minimal speed
$\displaystyle{c^{*}_{\Omega,A,q,f}(e)}$ is given by
\begin{equation}\label{min max ZFK }
\displaystyle{c^{*}_{\Omega,A,q,f}(e)}=\displaystyle{\min_{\varphi\in
E}\sup_{(s,x,y)\in\mathbb{R}\times\overline{\Omega}}R\,\varphi(s,x,y).}
\end{equation}
Furthermore, the $\min$ is attained by the function
$\phi^{*}(s,x,y)=u^{*}\left(\frac{s-x\cdot e}{c^{*}(e)},x,y\right)$ and its shifts $\phi^{*}(s+\tau,x,y)$ for any
$\tau\in\R,$ where $u^{*}$ is any solution of (\ref{front}) propagating with
the speed $c^{*}(e)=\displaystyle{c^{*}_{\Omega,A,q,f}(e)}.$
\end{thm}

In particular, Theorem \ref{thm min max ZFK} yields that formula (\ref{min max ZFK }) holds in the ``KPP'' case (\ref{KPP condition}) as well.

\begin{remark}
In Theorem \ref{thm min max max min comb}, the $\min$ and the $\max$ are attained by, and only by, the pulsating front $\phi(s,x,y)$ and its shifts
$\phi(s+\tau,x,y)$ for all $\tau\in\R.$ In Theorem \ref{thm min max ZFK}, the $\min$ is achieved by the front $\phi^{*}(s,x,y)$ with the speed
$c^{*}(e)$ and all its shifts $\phi^{*}(s+\tau,x,y).$ Actually, if the pulsating front $\phi^{*}$ is unique up to shift, then $\phi^{*}$ and its shifts are the unique minimizers in formula (\ref{min max ZFK }). The uniqueness is known in the ``KPP'' case (see Hamel, Roques \cite{Hamel Roques}), but it is still open in the general ``ZFK'' case.
\end{remark}

 We mention that a $\max$-$\min$ formula of the type ({\ref{max min
combustion}}) can not hold for the minimal speed $c^{*}(e)$ in the
``ZFK'' or the ``KPP'' case. A simple justification is given in section \ref{change of variables section}.

 The variational formulations of the speeds of
propagation which are given in Theorems \ref{thm min max max min comb} and \ref{thm min max ZFK} are more
general than those in Hamel \cite{Hamel} and Heinze, Papanicolaou,
Stevens \cite{HPS}. In Theorems \ref{thm min max max min comb} and \ref{thm min max ZFK}, we consider
nonhomogeneous nonlinearities $f=f(x,y,u)$ and the domain $\Omega$
is in the most general periodic situation. However, in \cite{Hamel}, the domain was an infinite cylinder of $\mathbb{R}^{N}$
and the advection $q$ was in the form of shear flows. Moreover, in this paper, the
nonhomogeneous operator $\nabla\cdot(A\nabla u)$ replaces the
Laplace operator $\Delta u$ taken in \cite{Hamel}. On the other hand, in \cite{HPS}, the domain $\Omega$ was an infinite cylinder in $\mathbb{R}^{N}$
with a bounded cross section. Namely,
$\Omega=\mathbb{R}\times\omega\subset\mathbb{R}^{N}$ where the
cross section $\omega$ is a bounded domain in $\mathbb{R}^{N-1}.$ Moreover, the authors did not consider an advection field in \cite{HPS}. Finally, concerning the
 nonlinearities, they were depending
only on $u$ (i.e $f=f(u)$ and is satisfying either (\ref{combustion
hom}) or (\ref{ZFK homogenous})) in both of \cite{Hamel} and \cite{HPS}.

Besides the fact that we consider here a wider family of diffusion and reaction
coefficients, our assumptions are less strict than those supposed in
\cite{HPS} and \cite{3volpert}. Roughly speaking, the authors, in \cite{HPS} and \cite{3volpert}, assume a
stability condition on the pulsating travelling fronts. We mention
that such a stability condition is fulfilled in the homogenous setting; however, it has not been rigorously proved so far that this condition is satisfied in the heterogenous setting. Meanwhile, the assumptions of the present paper only involve the coefficients of the reaction-advection-diffusion
equation (\ref{front}), and they can then be checked easily.

Actually, in the ``KPP'' case, another ``simpler'' variational formula for the minimal speed $ \displaystyle{c^*(e)=c^{*}_{\Omega,A,q,f}(e)}$ is known. This known formula involves only the linearized nonlinearity $f$ at $u=0. $ Namely, it follows from \cite{BHN1} that
\begin{thm}[Berestycki, Hamel, Nadirashvili \cite{BHN1}]\label{varthm} Let $e$ be a fixed unit vector in $\mathbb{R}^d$ and let
$\tilde{e}=(e,0,\ldots,0)\in\mathbb{R}^N.$ Assume that $f$ is a ``KPP'' nonlinearity and that $\Omega,A$ and $q$
satisfy (\ref{comega}), (\ref{cA}) and (\ref{cq}) respectively. Then, the minimal speed $c^{*}(e)$ of
pulsating fronts solving (\ref{front}) and propagating in the
direction of $e$ is given by
\begin{equation}\label{var}
    \displaystyle{c^*(e)=c^{*}_{\Omega,A,q,f}(e)\,=\,\min_{\lambda>0}\frac{k(\lambda)}{\lambda}},
\end{equation}
where $\displaystyle{k(\lambda)=k_{\Omega,e,A,q,\zeta}(\lambda)}$ is
the principal eigenvalue of the operator
$\displaystyle{L_{\Omega,e,A,q,\zeta,\lambda}}$ which is defined by
\begin{eqnarray}\label{Leq}
\begin{array}{ll}
\displaystyle{L_{\Omega,e,A,q,\zeta,\lambda}\psi\,:=}&\,\displaystyle{\nabla\cdot(A\nabla\psi)\,+2\lambda\tilde{e}\cdot
A\nabla\psi\,+q\cdot\nabla\psi\,}\\
&\displaystyle{+[\lambda^2\tilde{e}A\tilde{e}+\lambda\nabla\cdot(A\tilde{e})+\lambda
q\cdot\tilde{e}+\,\zeta]\psi}
\end{array}
\end{eqnarray}
acting on the set $$\begin{array}{ll}
\widetilde{E_{\lambda}}=&\left\{\psi\,\in\,C^2(\overline{\Omega}),
\psi\;\hbox{is $L$-periodic with respect to $x$ and}\right.\vspace{3 pt}\\
&\left.\ds{\;\;\nu\cdot
A\nabla\psi\,=\,-\lambda(\nu \cdot A
\tilde{e})\psi\;\hbox{on}\;\partial{\Omega}}\right\}.
\end{array}$$
\end{thm}

In our last result, we prove that formula (\ref{min max ZFK }) implies formula (\ref{var}) in the ``KPP'' case, but under some additional assumptions on the
advection and the diffusion coefficients. This result gives an alternate proof of the formula (\ref{var}).
\begin{thm}\label{min max implies min} Let $e$ be a fixed unit vector in $\mathbb{R}^d$ and let
$\tilde{e}=(e,0,\ldots,0)\in\mathbb{R}^N.$ Assume that $f$ is a ``KPP'' nonlinearity and that $\Omega,A$ and $q$
satisfy (\ref{comega}), (\ref{cA}) and (\ref{cq}) respectively. Assume, furthermore, that
$\nu\cdot A\te=0~\hbox{ on }~\partial\Omega$ (in the case where $\partial \Omega\not=\emptyset $).
Then, formula (\ref{min max ZFK }) implies formula (\ref{var}).
 \end{thm}

\section{Main tools: change of variables and maximum principles}\label{change
of variables section}
In this section, we introduce some tools that will be used in
different places of this paper in order to prove the main
results.

Throughout this paper, $\tilde{e}$ will denote the vector in
$\mathbb{R}^{N}$ defined by
$$\tilde{e}=(e,0,\cdots,0)=(e^{1},\cdots,e^{d},0,\cdots,0),$$
where $e^{1},\cdots,e^{d}$ are the components of the vector $e.$

Our study is concerned with the model (\ref{front}).
Having a ``combustion'', a ``ZFK'', or a  ``KPP'' nonlinearity,
together with the assumptions (\ref{cA}) and (\ref{cq}), problem
(\ref{front}) has at least a classical solution $(c,u)$ such that $c>0$ and $u_t>0$ {\rm(}see
Theorems \ref{Berest Hamel combustion Thm} and \ref{Berest Hamel ZFK
Thm}{\rm)}. The function $u$ is globally
$C^{1,\mu}(\mathbb{R}\times\overline{\Omega})$ and  $C^{2,\mu}$ with respect
to $(x,y)$ variables (for every $\mu\in[0,1)$). It follows that $\nabla_{x,y}.(A\nabla u)\in
C(\mathbb{R}\times\overline{\Omega}).$ Having a unit
direction $e\in\mathbb{R}^{d},$ and having a bounded classical
solution $(c,u)$ of (\ref{front}) with $c=c(e)$ {\rm(}combustion
case{\rm)} or $c\geq c^{*}(e)$ {\rm(}ZFK or KPP case{\rm)}, we make the same change of
variables as Xin \cite{Xin2}. Namely, let $\phi=\phi(s,x,y)$ be the
function defined by
\begin{equation}\label{change of variables}
\phi(s,x,y)=u\left(\frac{s-x\cdot e}{c},x,y\right)~\hbox{for
all}~s\in\mathbb{R}~\hbox{and}~(x,y)\in\overline{\Omega}.
\end{equation}
Then, for all $(s,x,y)\in\mathbb{R}\times\overline{\Omega},$
\begin{eqnarray*}
\begin{array}{cl}
\left[\nabla_{x,y}\cdot(A\nabla_{x,y}\phi)+(\tilde{e}A\tilde{e})\phi_{ss}+\nabla_{x,y}\cdot(A\tilde{e}\phi_s)\right.+\partial_s(\tilde{e}A\nabla_{x,y}\phi)\left]\right.(s,x,y)&\vspace{5 pt}\\
=\nabla_{x,y}\cdot(A\nabla u)(t,x,y),&
\end{array}
\end{eqnarray*}
where $s=x\cdot e+ct.$ Consequently,
\begin{eqnarray*}
F[\phi](s,x,y)=&\nabla_{x,y}\cdot(A\nabla_{x,y}\phi)+(\tilde{e}A\tilde{e})\phi_{ss}
+\nabla_{x,y}\cdot(A\tilde{e}\phi_s)
+\partial_s(\tilde{e}A\nabla_{x,y}\phi)
\end{eqnarray*}
is defined at each point
$(s,x,y)\in\mathbb{R}\times\overline{\Omega}$ and
 the map $(s,x,y)\mapsto F[\phi](s,x,y)$
belongs to $C(\mathbb{R}\times\overline{\Omega}).$

In all this paper, $L=L_{c}$ will denote the operator acting on the
set $E$ (given in Theorem \ref{thm min max max min comb}) and which is defined by
\begin{eqnarray}\label{degen elliptic eq sat by phi}
\begin{array}{ll}
L\varphi&=\nabla_{x,y}\cdot(A\nabla_{x,y}\varphi)+(\tilde{e}A\tilde{e})\varphi_{ss}+\nabla_{x,y}\cdot(A\tilde{e}\varphi_s)+\partial_s(\tilde{e}A\nabla_{x,y}\varphi)\vspace{4 pt}\\
&\,+\,q\cdot\nabla_{x,y}\varphi\,+\,(q\cdot\tilde{e}\,-\,c)\varphi_s~~\hbox{in}~C(\mathbb{R}\times\overline
{\Omega})\vspace{4 pt}\\
&=F[\varphi]+q\cdot\nabla_{x,y}\varphi+(q\cdot\tilde{e}\,-\,c)\varphi_s~~\hbox{in}~C(\R\times\overline{\Omega}),
\end{array}
\end{eqnarray}
for all $\varphi\in E.$

 It follows from above that if $\phi=\phi(s,x,y)$ is a
function that is given by a pulsating travelling $(c,u)$ solving
(\ref{front}) {\rm(}under the change of variables (\ref{change of
variables}){\rm)}, then $F[\phi]\in
C(\mathbb{R}\times\overline{\Omega}),$ $\phi$ is globally bounded in
$C^{1,\mu}(\mathbb{R}\times\overline{\Omega})$ (for every $\mu\in[0,1)$) and it satisfies the
following degenerate elliptic equation
\begin{eqnarray}\label{front reflected on phi}
\begin{array}{lll}
L\phi(s,x,y)+f(x,y,\phi)&=F[\phi](s,x,y)+\,q\cdot\nabla_{x,y}\phi(s,x,y)\,&\vspace{4 pt}\\
&+\,(q\cdot\tilde{e}\,-\,c)\phi_s(s,x,y)+f(x,y,\phi)&=0
\end{array}
\end{eqnarray}
in $\mathbb{R}\times\overline{\Omega},$ together with the boundary and periodicity conditions
\begin{eqnarray}\label{boundary and periodicity conditions on phi}
\left\{
  \begin{array}{ll}
    \phi~\hbox{is} \hbox{ $L-$periodic with respect
to $x,$}\\
\nu\cdot
A(\nabla_{x,y}\phi+\tilde{e}\phi_{s})=0~\hbox{on}~\mathbb{R}\times\overline{\Omega}.
  \end{array}
\right.
\end{eqnarray}
Moreover, since $u(t,x,y)\rightarrow0$ as
$x\cdot e\rightarrow-\infty$ and $u(t,x,y)\rightarrow1$ as $x\cdot
e\rightarrow+\infty$ locally in $t$ and uniformly in $y$ and in the directions of
$\mathbb{R}^{d}$ which are orthogonal to $e,$ and since $\phi$ is
$L-$periodic with respect to $x,$ the change of variables $s=x\cdot
e+ct$ guarantees that

\begin{equation}\label{phi(-infty,..)=0 phi(+infty,,..)=1}
\phi(-\infty,.,.)=0~\hbox{and}~\phi(+\infty,.,.)=1~\hbox{uniformly
in}~(x,y)\in\overline{\Omega}.
\end{equation}
Therefore, one can conclude that $\phi\in E.$

\begin{remark}
It is now clear that a max-min formula of the type (\ref{max min
combustion}) can not hold for the minimal speed $c^{*}(e)>0$ in the
``ZFK'' or the ``KPP'' case. Indeed, for each speed $c\geq
c^{*}(e),$ there is a solution $(c,u)$ of (\ref{front}) such that $u_t>0,$ which gives
birth to a function $\phi=\phi(s,x,y)$ under the change of variables
(\ref{change of variables}). Owing to the above discussions the function $\phi\in E$ and it satisfies
$$c=R\phi(s,x,y)~~\hbox{for
all}~~(s,x,y)\in\mathbb{R}\times\overline{\Omega}.$$ Therefore
$$\displaystyle{\sup_{\varphi\in
E}\inf_{(s,x,y)\in\mathbb{R}\times\overline{\Omega}}R\,\varphi(s,x,y)\geq
c}.$$ Since one can choose any $c\geq c^{*}(e),$ one concludes that
$$\displaystyle{\sup_{\varphi\in
E}\inf_{(s,x,y)\in\mathbb{R}\times\overline{\Omega}}R\,\varphi(s,x,y)=+\infty}$$
in the ``ZFK'' or the ``KPP'' case.
\end{remark}

\begin{remark}[The same formul{\ae}, but over a subset of $E$]
If the restriction of the nonlinear source $f$ in (\ref{front}) is
$C^{1,\delta}(\overline{\Omega}\times[0,1]),$ one can then conclude
that (see the proof of Proposition 6.3 in \cite{BH1}) any solution
$u$ of (\ref{front}) satisfies:
$$\forall(t,x,y)\in\mathbb{R}\times\overline{\Omega},~~|\partial_{tt}u(t,x,y)|\leq M\,\partial_{t}u(t,x,y)$$
for some constant $M$ independent of $(t,x,y).$ In other words, the
function $$\phi(s,x,y)=u((s-x\cdot e)/c,x,y)$$ (where $c=c(e)$ in
the ``combustion'' case, and $c=c^{*}(e)$ in the ``ZFK'' or the
``KPP'' case) satisfies
$$\forall(s,x,y)\in\mathbb{R}\times\overline{\Omega},~~|\partial_{ss}\phi(s,x,y)|\leq
({M}/{c})\;\partial_{s}\phi(s,x,y).$$ Let $E^{\,'}$ be the
functional subset of $E$ defined by
$$E^{\,'}=\left\{\varphi\in~E,~\exists C>0,~\forall \,(s,x,y)\in\mathbb{R}\times\overline{\Omega},\,|\partial_{ss}\varphi(s,x,y)|\leq
C\;\partial_{s}\varphi(s,x,y)\right\}.$$ The
previous facts together with the discussions at the beginning of this section
imply that the functions $\phi$ and $\phi^{*}$ of Theorems \ref{thm min max max min comb} and \ref{thm min max ZFK} are elements of $E^{\,'}\subset E.$ These theorems also
yield that the $\max$-$\min$ and the $\min$-$\max$ formul{\ae} can
also hold over the subset $E^{\,'}$ of $E.$

Namely, in the case of a ``combustion'' nonlinearity
\begin{equation}\label{min max comb on E'}
c(e)= \displaystyle{\min_{\varphi\in
E^{\,'}}\sup_{(s,x,y)\in\mathbb{R}\times\overline{\Omega}}R\,\varphi(s,x,y)}
\end{equation}
and
\begin{equation}\label{max min comb on E'}
 c(e)= \displaystyle{\max_{\varphi\in
E^{\,'}}\inf_{(s,x,y)\in\mathbb{R}\times\overline{\Omega}}R\,\varphi(s,x,y)}.
\end{equation}
Moreover, the $\min$ and the $\max$ are attained at, and only at, the function
$\phi(s,x,y)$ and its shifts $\phi(s+\tau,x,y)$ for any $\tau\in\R.$

 On the other hand, only a $\min$-$\max$ formula holds in the
case of ``ZFK'' or ``KPP'' nonlinearities. That is
\begin{equation}\label{min max ZFK over E'}
c^{*}(e)= \displaystyle{\min_{\varphi\in
E^{\,'}}\sup_{(s,x,y)\in\mathbb{R}\times\overline{\Omega}}R\,\varphi(s,x,y)}.
\end{equation}
Moreover, the $\min$ is attained at the function $\phi^{*}(s,x,y)$ and its shifts $\phi^{*}(s+\tau,x,y)$ for any
$\tau\in\R.$
\end{remark}

 In the proofs of the variational formul{\ae} which were given in
Theorem \ref{thm min max max min comb} and Theorem \ref{thm min max ZFK}, we will use two versions of the maximum
principle in unbounded domains for some problems related to
(\ref{degen elliptic eq sat by phi}-\ref{boundary and periodicity
conditions on phi}) and (\ref{phi(-infty,..)=0 phi(+infty,,..)=1}).
Such generalized maximum principles were proved in Berestycki, Hamel
\cite{BH1} in a slightly more general framework:

\begin{lemm}[\cite{BH1}]\label{max principle1}
Let $e$ be a fixed unit vector in $\mathbb{R}^{d}.$ Let $g(x,y,u)$
be a globally bounded and globally Lipschitz-continuous function
defined in $\overline{\Omega}\times\mathbb{R}$ and assume that $g$
is non-increasing with respect to $u$ in
$\overline{\Omega}\times(-\infty,\delta]$ for some $\delta>0.$ Let
$h\in\mathbb{R}$ and
$\displaystyle{\Sigma^{-}_{h}:=(-\infty,h)\times\Omega.}$ Let
$c\neq0$ and $\phi^{1}(s,x,y),\;\phi^{2}(s,x,y)$ be two bounded and
globally
$\displaystyle{C^{1,\mu}\left(\overline{\Sigma^{-}_{h}}\right)}$
functions (for some $\mu>0$) such that
\begin{eqnarray}\label{L phi1 >0, L phi2 <0}
\left\{
  \begin{array}{rcl}
L\,\phi^{1}+g(x,y,\phi^{1})&\geq&0~\hbox{ in }~\mathcal{D}{'}(\Sigma^{-}_{h}),\vspace{3 pt}\\
L\,\phi^{2}+g(x,y,\phi^{2})&\leq&0~\hbox{ in }~\mathcal{D}{'}(\Sigma^{-}_{h}),\vspace{3 pt}\\
\nu\cdot
A\left[\tilde{e}(\phi^{1}_{s}-\phi^{2}_{s})+\nabla_{x,y}(\phi^{1}-\phi^{2})\right]&\leq&0\hbox{ on }(-\infty,h]\times\partial\Omega,\vspace{3 pt}\\
\displaystyle{\lim_{{s_0}\rightarrow-\infty}\sup_{\{s\leq
s_{0},\;(x,y)\in\overline{\Omega}\}}\,[\phi^{1}(s,x,y)-\phi^{2}(s,x,y)]}&\leq&0,
\end{array}
\right.
\end{eqnarray}
where \begin{eqnarray}\label{the operator L of the change ofvariables}
 \begin{array}{lll}
L\,\phi&:=&\,\nabla_{x,y}\cdot(A\nabla_{x,y}\phi)+(\tilde{e}A\tilde{e})\phi_{ss}+\nabla_{x,y}\cdot(A\tilde{e}\phi_s)+\partial_s(\tilde{e}A\nabla_{x,y}\phi)\\
&&\,+\,q\cdot\nabla_{x,y}\phi\,+\,(q\cdot\tilde{e}\,-\,c)\phi_s,
\end{array}
\end{eqnarray}
and $\tilde{e}$ denotes the vector
$(e,0,\cdots,0)\in\mathbb{R}^{N}.$

If $\phi^{1}\leq\delta$ in $\overline{\Sigma^{-}_{h}}$ and
$\phi^{1}(h,x,y)\leq\phi^{2}(h,x,y)$ for all
$(x,y)\in\overline{\Omega},$ then
$$\phi^{1}\leq\phi^{2}\quad\hbox{in}\;\overline{\Sigma^{-}_{h}}.$$
\end{lemm}
\begin{remark}
Note here that $\phi^{1},\,\phi^{2},\,q,\,A$ and $g$ are not assumed
to be $L-$periodic in $x$ and that $q$ is not assumed to satisfy
(\ref{cq}).
\end{remark}

Changing $\phi^{1}(s,x,y),\,\phi^{2}(s,x,y)$ and $g(x,y,s)$ into $1-\phi^{1}(-s,x,y),\,1-\phi^{2}(-s,x,y)$ and $-g(x,y,1-s)$ respectively in Lemma \ref{max principle1} leads to the
following

\begin{lemm}[\cite{BH1}]\label{maxprinciple2}
Let $e$ be a fixed unit vector in $\mathbb{R}^{d}.$ Let $g(x,y,u)$
be a globally bounded and globally Lipschitz-continuous function
defined in $\overline{\Omega}\times\mathbb{R}$ and assume that $g$
is non-increasing with respect to $u$ in
$\overline{\Omega}\times[1-\delta,+\infty)$ for some $\delta>0.$ Let
$h\in\mathbb{R}$ and
$\displaystyle{\Sigma^{+}_{h}:=(h,+\infty)\times\Omega.}$ Let
$c\neq0$ and $\phi^{1}(s,x,y),\;\phi^{2}(s,x,y)$ be two bounded and
globally
$\displaystyle{C^{1,\mu}\left(\overline{\Sigma^{+}_{h}}\right)}$
functions (for some $\mu>0$) such that
\begin{eqnarray}\label{L phi1 >0, L phi2 <0}
\left\{
  \begin{array}{rcl}
L\,\phi^{1}+g(x,y,\phi^{1})&\geq&0~\hbox{ in }~\mathcal{D}{'}(\Sigma^{+}_{h}),\vspace{3 pt}\\
L\,\phi^{2}+g(x,y,\phi^{2})&\leq&0~\hbox{ in }~\mathcal{D}{'}(\Sigma^{+}_{h}),\vspace{3 pt}\\
\nu\cdot
A\left[\tilde{e}(\phi^{1}_{s}-\phi^{2}_{s})+\nabla_{x,y}(\phi^{1}-\phi^{2})\right]&\leq&0\hbox{ on }[h,+\infty)\times\partial\Omega,\vspace{3 pt}\\
\displaystyle{\lim_{{s_0}\rightarrow+\infty}\sup_{\{s\geq
s_{0},\;(x,y)\in\overline{\Omega}\}}\,[\phi^{1}(s,x,y)-\phi^{2}(s,x,y)]}&\leq&0,
\end{array}
\right.
\end{eqnarray}
where $L$ is the same operator as in Lemma \ref{max principle1}.

If $\phi^{2}\geq 1-\delta$ in $\overline{\Sigma^{+}_{h}}$ and
$\phi^{1}(h,x,y)\leq\phi^{2}(h,x,y)$ for all
$(x,y)\in\overline{\Omega},$ then
$$\phi^{1}\leq\phi^{2}\quad\hbox{in}~~\overline{\Sigma^{+}_{h}}.$$
\end{lemm}

\section{ Case of a ``combustion'' nonlinearity}
This section is devoted to prove Theorem \ref{thm min max max min comb}, where the nonlinearity $f$ satisfies the assumptions
(\ref{f1}-\ref{f2}) and (\ref{combustion}).

\subsection{ Proof of formula (\ref{min max combustion})}\label{proof of min max combustion}

 Having a prefixed unit
direction $e\in\mathbb{R}^{d},$ and since the coefficients $A$ and
$q$ of problem (\ref{front}) satisfy the assumptions (\ref{cA}) and
(\ref{cq}), it follows, from Theorem \ref{Berest Hamel combustion
Thm}, that there exists a unique pulsating travelling front
$(c(e),u)$ ($u$ is unique up to a translation in the time variable)
which solves problem (\ref{front}). Moreover, $\partial_{t}u>0$ in
$\mathbb{R}\times\overline{\Omega}.$ We will complete the proof of (\ref{min max combustion}) via
two steps.

\textit{Step 1.} After the discussions done in the section
\ref{change of variables section}, the existence of a classical
solution $(c(e),u),$ satisfying (\ref{front}), implies the
existence of a globally $C^{1}(\mathbb{R}\times\overline{\Omega})$
function $\phi(s,x,y)$ satisfying $0\leq\phi\leq1$ in
$\mathbb{R}\times\overline{\Omega},$ with
\begin{eqnarray}\label{eqs sat by phi}
\left\{
  \begin{array}{ll}
\phi~\hbox{is $L-$periodic with respect to~} x,\\
 L\phi(s,x,y)+f(x,y,\phi)=0~\hbox{in}~\mathcal{D}'(\mathbb{R}\times\overline{\Omega}),\\
\nu\cdot
A(\nabla_{x,y}\phi+\tilde{e}\phi_{s})=0~\hbox{in}~\mathbb{R}\times\partial\Omega,\\
\phi(-\infty,.,.)=0, ~\hbox{and}~\phi(+\infty,.,.)=1~\hbox{uniformly
in}~(x,y)\in\overline{\Omega},
  \end{array}
\right.
\end{eqnarray}
where $L$ is the operator defined in (\ref{degen elliptic eq sat by
phi}) for $c=c(e).$ We also recall that the two functions $u$ and $\phi$ satisfy the
relation $$u(t,x,y)=\phi(x\cdot
e+c(e)t,x,y),~~(t,x,y)\in\mathbb{R}\times\overline{\Omega}.$$
One has $\partial_{s}\phi>0$ in
$\mathbb{R}\times\overline{\Omega}$ and this is equivalent to say
that the function $u=u(t,x,y)$ is increasing in $t,$ since $c(e)>0.$

Together with the facts in section \ref{change of variables}, one gets that
the function $\phi\in E.$ Furthermore, (\ref{eqs sat by phi}) yields
that
\begin{equation}\label{c(e)=Rphi+q.e}
\forall\,s\in\mathbb{R},\;\forall(x,y)\in\overline{\Omega},\;c(e)=R\,\phi(s,x,y),
\end{equation}
and
\begin{equation}\label{Lphi+f(x,yphi)=0}
    L\phi(s,x,y)+f(x,y,\phi)=0,
\end{equation}
 where $R\phi$ is the function defined in Theorem \ref{thm min max max min comb}. In other words, the $L-$periodic (with respect to $x$)
function $R\phi\;$ is constant over
$\mathbb{R}\times\overline{\Omega}$ and it is equal to $c(e).$

It follows, from (\ref{c(e)=Rphi+q.e}) and from the above
explanations, that
$$c(e)\geq\inf_{\varphi\in E}\sup_{(s,x,y)\in\mathbb{R}\times\overline{\Omega}}R\varphi(s,x,y).$$

To complete the proof of formula (\ref{min max combustion}), we
assume that $$c(e)>\inf_{\varphi\in
E}\sup_{(s,x,y)\in\mathbb{R}\times\overline{\Omega}}R\varphi(s,x,y).$$
Then, there exists a function $\psi=\psi(s,x,y)\in E$ such that
$$c(e)>\sup_{(s,x,y)\in\mathbb{R}\times\overline{\Omega}} R\psi\,(s,x,y).$$
Since the function $\psi\in E,$ one then has $\psi_{s}(s,x,y)>0$ for
all $(s,x,y)\in\mathbb{R}\times\overline{\Omega}.$ This yields that
\begin{equation}\label{Lpsi+...<0}
L\psi(s,x,y)+f(x,y,\psi)<0~\hbox{in}~\mathbb{R}\times\overline{\Omega},
\end{equation}
where $L$ is the operator defined in (\ref{degen elliptic eq sat by
phi}) for $c=c(e).$

Notice that the later holds for each function of the type
$$\psi^{\tau}(s,x,y):=\psi(s+\tau,x,y)$$ because of the invariance of
(\ref{Lpsi+...<0}) with respect to $s$ and because the advection
field $q$ and the diffusion matrix $A$ depend on the variables
$(x,y)$ only. That is
\begin{equation}\label{Lpsi^tau}
L\psi^{\tau}(s,x,y)+f(x,y,\psi^{\tau})<0~\hbox{in}~\mathbb{R}\times\overline{\Omega}.
\end{equation}

\textit{Step 2.} In order to draw a contradiction, we are going to
slide the function $\psi$ with respect to $\phi.$ From the limiting
conditions satisfied by these two functions, there exists a real
number $B>0$ such that
\begin{eqnarray*}
\left\{
  \begin{array}{ll}
    \phi(s,x,y)\leq\theta&\hbox{for
all}~s\leq-B,~(x,y)\in\overline{\Omega},\\
\psi(s,x,y)\geq1-\rho&\hbox{for all}~s\geq
B,~(x,y)\in\overline{\Omega},
  \end{array}
\right.
\end{eqnarray*}
and
\begin{equation}\label{phi(B,...)}
\phi(B,x,y)\geq1-\rho~~\hbox{for all}~(x,y)\in\overline{\Omega},
\end{equation}
where $\theta$ and $\rho$ are the values that appear in the
conditions (\ref{combustion}) satisfied by the ``combustion''
nonlinearity $f.$ Taking $\tau\geq 2B,$ and since $\psi$ is
increasing with respect to $s,$ one gets that
$\phi(-B,x,y)\leq\psi^{\tau}(-B,x,y)$ for all
$(x,y)\in\overline{\Omega}$ and $\psi^{\tau}\geq 1-\rho$ in
$\overline{\Sigma_{-B}^{+}}.$

 It follows from Lemma \ref{max principle1} (take $\delta=\theta,~h=-B,~\phi^{1}=\phi,~\hbox{and}~\phi^{2}=\psi^{\tau}$) that
$\phi\leq\psi^{\tau}$ in $\overline{\Sigma_{-B}^{-}}.$ Moreover,
Lemma \ref{maxprinciple2} (take
$\delta=\rho,~h=-B,~\phi^{1}=\phi,~\hbox{and}~\phi^{2}=\psi^{\tau}$
) implies that $\phi\leq\psi^{\tau}$ in
$\overline{\Sigma_{-B}^{+}}.$ Consequently, $\phi\leq\psi^{\tau}$ in
$\mathbb{R}\times\overline{\Omega}\,$ for all $\tau\geq 2B.$

Let us now decrease $\tau$ and set
$$\tau^{*}=\inf\{\tau\in\mathbb{R},~~\phi\leq\psi^{\tau}~\hbox{in}~\mathbb{R}\times\overline{\Omega}~\}.$$
 First one notes that $\tau^{*}\leq2B.$ On the other hand, the limiting
conditions $\psi(-\infty,.,.)=0$ and $\phi(+\infty,.,.)=1$
imply that $\tau^{*}$ is finite. By continuity,
$\phi\leq\psi^{\tau^*}~\hbox{in}~\mathbb{R}\times\overline{\Omega}.$
Two cases may occur according to the value of $\displaystyle{\sup_{[-B,B]\times\overline{\Omega}}\left(\phi-\psi^{\tau^*}\right).}$

\underline{\textit{case 1:}} suppose that
$$\displaystyle{\sup_{[-B,B]\times\overline{\Omega}}\left(\phi-\psi^{\tau^*}\right)<0}.$$
Since the functions $\psi$ and $\phi$ are globally
$C^{1}(\mathbb{R}\times\overline{\Omega})$ there exists $\eta>0$ such that the above inequality
holds for all $\tau\in[\tau^*-\eta,\tau^*].$ Choosing any $\tau$ in
the interval $[\tau^*-\eta,\tau^*],$ and applying Lemma \ref{max
principle1} to the functions $\psi^{\tau}$ and $\phi,$ one gets that
$$\phi(s,x,y)\leq\psi^{\tau}(s,x,y)~~\hbox{for
all}~s\leq-B,~(x,y)\in\overline{\Omega},$$ together with the
inequality $$\phi(s,x,y)<\psi^{\tau}(s,x,y)~~\hbox{for
all}~s\in[-B,B],~\hbox{and for all}~(x,y)\in\overline{\Omega.}$$
Owing to (\ref{phi(B,...)}) and to the above inequality, it follows
that $$\psi^{\tau}(B,x,y)\geq 1-\rho~\hbox{in}~\overline{\Omega}.$$
Moreover, since the function $\psi$ is increasing in $s,$ one gets
that $\psi^{\tau}\geq1-\rho$ in $\overline{\Sigma_{B}^{+}}.$ Lemma
\ref{maxprinciple2}, applied to $\phi$ and $\psi^{\tau}$ in
$\overline{\Sigma_{B}^{+}},$ yields that
$$\phi(s,x,y)\leq\psi^{\tau}(s,x,y)~\hbox{for all}~s\geq B,~(x,y)\in\overline{\Omega}.$$
As a consequence, one obtains $\phi\leq\psi^{\tau}$ in
$\mathbb{R}\times\overline{\Omega},$ and that contradicts the
minimality of $\tau^{*}.$ Therefore, case 1 is ruled out.

\underline{\textit{case 2:}} suppose that
$$\displaystyle{\sup_{[-B,B]\times\overline{\Omega}}\left(\phi-\psi^{\tau^*}\right)=0}.$$
Then, there exists a sequence of points $(s_n,x_n,y_n)$ in
$[-B,B]\times\overline{\Omega}$ such that
$$\phi(s_n,x_n,y_n)-\psi^{\tau}(s_n,x_n,y_n)\rightarrow 0~\hbox{as}~n\rightarrow+\infty.$$
Due to the $L-$ periodicity of the functions $\phi$ and $\psi,$ one
can assume that $(x_n,y_n)\in\overline{C}.$ Consequently, one can
assume, up to extraction of a subsequence, that
$(s_n,x_n,y_n)\rightarrow(\bar{s},\bar{x},\bar{y})\in[-B,B]\times\overline{C}$
as $n\rightarrow+\infty.$ By continuity, one gets
$\phi(\bar{s},\bar{x},\bar{y})=\psi^{\tau^*}(\bar{s},\bar{x},\bar{y}).$

We return now to the variables $(t,x,y).$ Let
\begin{eqnarray*}
  z(t,x,y)&=&\phi(x\cdot e+c(e)\,t,x,y)-\psi(x\cdot
e+c(e)\,t+\tau^{*},x,y)\\
&=&u(t,x,y)-\psi(x\cdot e+c(e)\,t+\tau^{*},x,y) \hbox{ for all } (t,x,y)\in\mathbb{R}\times\overline{\Omega}.
\end{eqnarray*}
 Since the functions $\phi$ and $\psi$ are in $E,$ it follows that
the function $z$ is globally
$C^{1}(\mathbb{R}\times\overline{\Omega})$ and it satisfies
$$\forall\,(t,x,y)\in\mathbb{R}\times\Omega,~~\nabla_{x,y}\cdot(A\nabla z)(t,x,y)=F[\phi](s,x,y)-F[\psi^{\tau^{*}}](s,x,y),$$
where $s=x\cdot e+c(e)t.$ Thus, $\nabla_{x,y}\cdot(A\nabla z)\in
C(\mathbb{R}\times\overline{\Omega}).$ Moreover, the function $z$ is
non positive and it vanishes at the point $((\bar{s}-\bar{x}\cdot
e)/c(e),\bar{x},\bar{y}).$ It satisfies the boundary condition
$\nu\cdot(A\nabla z)=0$ on $\mathbb{R}\times\partial\Omega.$
Furthermore, it follows, from (\ref{c(e)=Rphi+q.e}) and
(\ref{Lpsi+...<0}), that
$$\partial_{t}z-\nabla_{x,y}\cdot (A\nabla z)+q(x,y)\cdot\nabla_{x,y}z\leq f(x,y,\phi)-f(x,y,\psi^{\tau^{*}}).$$
However, the function $f$ is globally Lipschitz-continuous in
$\overline{\Omega}\times\mathbb{R};$ hence, there exists a bounded
function $b(t,x,y)$ such that $$\partial_{t}z-\nabla_{x,y}\cdot
(A\nabla z)+q(x,y)\cdot\nabla_{x,y}z+b(t,x,y)\,z\leq
0~\hbox{in}~\mathbb{R}\times\Omega,$$ with $z(t,x,y)\leq0$ for all $(t,x,y)\in\R\times\overline{\Omega}.$

Applying the strong parabolic maximum principle and Hopf lemma, one
gets that $z(t,x,y)=0$ for all $t\leq (\bar{s}-\bar{x}\cdot e)/c(e)$
and for all $(x,y)\in \overline{\Omega}.$ On the other hand, it
follows from the definition of $z$ and from the $L-$periodicity of
the functions $\phi$ and $\psi$ that $z(t,x,y)=0$ for all
$(t,x,y)\in\mathbb{R}\times\overline{\Omega}.$ Consequently,
$$\phi(s,x,y)=\psi^{\tau^*}(s,x,y)=\psi(s+\tau^*,x,y)~\hbox{for all}~(s,x,y)\in\mathbb{R}\times\overline{\Omega}.$$

Referring to the equations (\ref{Lphi+f(x,yphi)=0}) and
(\ref{Lpsi^tau}), one gets a contradiction. Thus, case 2 is
ruled out too, and that completes the proof of the formula (\ref{min max
combustion}).
\begin{remark}[The uniqueness, up to a shift, of the minimizer in (\ref{min max combustion})]\label{minimizer is unique up to shift}
If $\psi\in E$ is a minimizer in (\ref{min max combustion}). The above arguments imply that case 2 necessarily occurs, and that $\psi$ is equal to a shift of $\phi.$ In other words, the minimum in (\ref{min max combustion}) is realized by and only by the shifts of $\phi.$
\end{remark}

\subsection{ Proof of formula (\ref{max min combustion})}
In this subsection, we are going to prove the ``max-min'' formula of
the speed of propagation $c(e)$ whenever the nonlinearity $f$ is of
the ``combustion'' type. The tools and techniques which one uses
here are similar to those used in the previous subsection. However,
we are going to sketch the proof of formula (\ref{max min
combustion}) for the sake of completeness.

As it was justified in the previous subsection, one easily gets that
$$c(e)\leq\sup_{\varphi\in E}\inf_{(s,x,y)\in\mathbb{R}\times\overline{\Omega}}R\varphi(s,x,y)$$
and
$$\forall\,(s,x,y)\in\mathbb{R}\times\overline{\Omega},\quad c(e)=R\phi(s,x,y),$$
where $$\phi(s,x,y)=u\left((s-x\cdot e)/c(e),x,y\right),~\hbox{for
all}~(s,x,y)\in\mathbb{R}\times\overline{\Omega},$$  and
$u=u(t,x,y)$ is the unique (up to a translation in $t$) pulsating
travelling front solving problem (\ref{front}) and propagating in the
speed $c(e).$ We recall that the function $\phi\in E$ (see
section \ref{change of variables section}). It follows that the function $\phi$
satisfies the following
\begin{eqnarray}\label{eqs sat by phi 1}
\left\{
  \begin{array}{ll}
\phi~\hbox{is $L-$periodic with respect to~} x,\\
\phi~\hbox{is increasing in~} s\in\mathbb{R},\\
 L\phi(s,x,y)+f(x,y,\phi)=0~\hbox{in}~\mathbb{R}\times\overline{\Omega},\\
\nu\cdot
A(\nabla_{x,y}\phi+\tilde{e}\phi_{s})=0~\hbox{in}~\mathbb{R}\times\partial\Omega,\\
\phi(-\infty,.,.)=0, ~\hbox{and}~\phi(+\infty,.,.)=1~\hbox{uniformly
in}~(x,y)\in\overline{\Omega},
  \end{array}
\right.
\end{eqnarray}
where $L$ is the operator defined in (\ref{degen elliptic eq sat by
phi}) for $c=c(e).$

Notice that the later holds also for each function of the type
$$\phi^{\tau}(s,x,y):=\phi(s+\tau,x,y)$$ because of the invariance of
(\ref{Lpsi+...>0}) with respect to $s$ and because the advection
field $q$ and the diffusion matrix $A$ depend on the variables
$(x,y)$ only.

To complete the proof of formula (\ref{max min combustion}), we
assume that $$c(e)<\sup_{\varphi\in
E}\inf_{(s,x,y)\in\mathbb{R}\times\overline{\Omega}}R\varphi(s,x,y).$$
Hence, there exists $\psi\in E$ such that
$$c(e)< R\psi\,(s,x,y),~\hbox{for all}~(s,x,y)\in\mathbb{R}\times\overline{\Omega}.$$
Since the function $\psi\in E,$ one then has $\psi_{s}(s,x,y)>0$ for
all $(s,x,y)\in\mathbb{R}\times\overline{\Omega}.$ This yields that
\begin{equation}\label{Lpsi+...>0}
L\psi(s,x,y)+f(x,y,\psi)>0~\hbox{in}~\mathbb{R}\times\overline{\Omega}.
\end{equation}

To get a contradiction, we are going to slide the function $\phi$
with respect to $\psi.$ In fact, the limiting conditions satisfied
by $\psi$ and $\phi,$ which are elements of $E,$ yield that there
exists a real positive number $B$ such that

\begin{eqnarray*}
\left\{
  \begin{array}{ll}
    \psi(s,x,y)\leq\theta&\hbox{for
all}~s\leq-B,~(x,y)\in\overline{\Omega},\\
\phi(s,x,y)\geq1-\rho&\hbox{for all}~s\geq
B,~(x,y)\in\overline{\Omega},
  \end{array}
\right.
\end{eqnarray*}
and
\begin{equation}\label{psi(B,...)}
\psi(B,x,y)\geq1-\rho~~\hbox{for all}~(x,y)\in\overline{\Omega},
\end{equation}
where $\theta$ and $\rho$ are the values appearing in the conditions
(\ref{combustion}) satisfied by the nonlinearity $f.$ Having
$\tau\geq2B,$ one applies Lemma \ref{max principle1} (taking
$\delta=\theta,~h=-B,~\phi^{1}=\psi,~\hbox{and}~\phi^{2}=\phi^{\tau}$)
and Lemma \ref{maxprinciple2} (taking
$\delta=\rho,~h=-B,~\phi^{1}=\psi,~\hbox{and}~\phi^{2}=\phi^{\tau}$)
to the functions $\phi^{\tau}$ and $\psi,$ over the domains
$\Sigma_{-B}^{-}$ and $\Sigma_{-B}^{+}$ respectively, to get that
$\psi\leq\phi^{\tau}$ in $\Sigma_{-B}^{-}$ and $\psi\leq\phi^{\tau}$
in $\Sigma_{-B}^{+}.$ Consequently, one can conclude that$$\forall
\,\tau\geq2B,~~\psi\leq\phi^{\tau}~~\hbox{in}~\mathbb{R}\times\overline{\Omega}.$$

Let us now decrease $\tau$ and
set$$\tau^{*}=\inf\{\tau\in\mathbb{R},~~\psi\leq\phi^{\tau}~\hbox{in}~\mathbb{R}\times\overline{\Omega}~\}.$$
It follows, from the limiting conditions $\psi(+\infty,.,.)=1$ and
$\phi(-\infty,.,.)=0,$ that $\tau^{*}$ is finite. By continuity, we
have $\psi\leq\phi^{\tau^*}.$ In this situation, two cases may
occur. Namely,

$$\hbox{case A:}\hskip2cm\displaystyle{\sup_{[-B,B]\times\overline{\Omega}}\left(\psi-\phi^{\tau^*}\right)<0},$$
or
$$\hbox{case B:}\hskip2cm\displaystyle{\sup_{[-B,B]\times\overline{\Omega}}\left(\psi-\phi^{\tau^*}\right)=0}.$$
Imitating the ideas and the skills used in case 1 and case 2 during
the proof of formula (\ref{min max combustion}), one gets that case
A (owing to minimality of $\tau^{*}$) and case B (owing to (\ref{eqs
sat by phi 1}) and (\ref{Lpsi+...>0})) are ruled out.

Therefore, the assumption that $$c(e)<\sup_{\varphi\in
E}\inf_{(s,x,y)\in\mathbb{R}\times\overline{\Omega}}R\varphi(s,x,y)$$
is false, and that completes the proof of formula (\ref{max min
combustion}).
\begin{remark}[The uniqueness, up to a shift, of the maximizer in (\ref{max min combustion})]
Similar to what we have already mentioned in Remark \ref{minimizer is unique up to shift}, if $\psi\in E$ is a maximizer in (\ref{max min combustion}), then the above arguments yield that case B necessarily occurs, and that $\psi$ is equal to a shift of $\phi.$ One then concludes that the maximum in (\ref{max min combustion}) is realized by, and only by, the shifts of $\phi.$
\end{remark}
\section{Case of  ``ZFK'' or  ``KPP'' nonlinearities: proof of formula (\ref{min max ZFK })} \label{section pf of min max zfk}
This section is devoted to the proof of Theorem \ref{thm min max ZFK}. We assume that the nonlinear source $f$ is of
``ZFK'' type. Remember that this case includes the class of ``KPP'' nonlinearities.
Namely, $f=f(x,y,u)$ is a nonlinearity satisfying
(\ref{f1}-\ref{f2}) and (\ref{ZFK}). We will divide the proof of formula (\ref{min max
ZFK }) into 3 steps:

\textit{Step 1.} Under the assumptions (\ref{comega}), (\ref{cA}),
and (\ref{cq}) on the domain $\Omega,$ the diffusion matrix $A,$ and
the advection field $q$ respectively, and having a
 nonlinearity $f$ satisfying the above assumptions, Theorem \ref{Berest Hamel ZFK Thm} yields that for
$c=\displaystyle{c^{*}_{\Omega,A,q,f}(e)},$ there exists a solution
$u^{*}=u^{*}(t,x,y)$ of (\ref{front}) such that $u^{*}_t(t,x,y)>0$
for all $(t,x,y)\in\mathbb{R}\times\Omega.$ In other words, the
function $\phi^{*}$ defined by
$$\phi^{*}(s,x,y)=u^{*}\left(\frac{s-x\cdot e}{c^{*}(e)},x,y\right),~(s,x,y)\in\mathbb{R}\times\overline{\Omega}$$
is increasing in $s\in\mathbb{R}.$ Owing to section \ref{change of variables section}, $\phi^{
*}$ satisfies
\begin{eqnarray}\label{degen elliptic eq sat by phi*}
\begin{array}{rl}
F[\phi^{*}]+\,q\cdot\nabla_{x,y}\phi^{*}\,+\,(q\cdot\tilde{e}\,-\,c^{*}(e))\phi^{*}_s,+f(x,y,\phi^{*})&=0~~\hbox{in}~~\mathbb{R}\times\overline{\Omega}
\end{array}
\end{eqnarray}
 together with boundary and periodicity conditions
\begin{eqnarray}\label{boundary and periodicity conditions on phi*}
\left\{
  \begin{array}{ll}
    \phi^{*}~\hbox{is} \hbox{ $L-$periodic with respect
to $x,$}\\
\nu\cdot
A(\nabla_{x,y}\phi^{*}+\tilde{e}\phi^{*}_{s})=0~\hbox{on}~\mathbb{R}\times\overline{\Omega}.
  \end{array}
\right.
\end{eqnarray}
 Moreover, (\ref{degen elliptic eq sat by phi*})
implies that
\begin{eqnarray}\label{Rphi*+q.e is constant}
\begin{array}{rl}
\forall~(s,x,y)\in\mathbb{R}\times\overline{\Omega},&\\
c^{*}(e)=&\displaystyle{\frac{F[\phi^{*}](s,x,y)+\,q\cdot\nabla_{x,y}\phi^{*}(s,x,y)+f(x,y,\phi^{*})}{\partial_{s}\phi^{*}(s,x,y)}+q(x,y)\cdot\tilde{e}}\\
=&R\phi^{*}(s,x,y),
\end{array}
\end{eqnarray}
and hence
$$c^{*}(e)\geq\displaystyle{\inf_{\varphi\in E}\sup_{(s,x,y)\in\mathbb{R}\times\overline{\Omega}}\frac{F[\varphi](s,x,y)+\,q\cdot\nabla_{x,y}\varphi+f(x,y,\varphi)}{\partial_{s}\phi(s,x,y)}+q(x,y)\cdot\tilde{e}}.$$

In order to prove equality, we argue by contradiction. Assuming that
the above inequality is strict, one can find $\delta>0$ such that
\begin{equation}\label{strict ineq}
c^{*}(e)-\delta>\displaystyle{\inf_{\varphi\in
E}\sup_{(s,x,y)\in\mathbb{R}\times\overline{\Omega}}\frac{F[\varphi](s,x,y)+\,q\cdot\nabla_{x,y}\varphi+f(x,y,\varphi)}{\partial_{s}\varphi(s,x,y)}+q(x,y)\cdot\tilde{e}}.
\end{equation}
To draw a contradiction, we are going to approach the ``ZFK''
nonlinearity $f$ by a sequence of ``combustion'' nonlinearities
$(f_\theta)_{\theta}$ and the minimal speed of propagation by the
sequence of speeds $(c_\theta)_{\theta}$ corresponding to the
functions $(f_\theta)_{\theta}.$ The details will appear in the next
step. \vskip0.4cm

\textit{Step 2.} Let $\chi$ be a $C^{1}(\mathbb{R})$ function such
that $0\leq\chi\leq1$ in $\mathbb{R},$ $\chi(u)=0$ for all $u\leq1,$
$0<\chi(u)<1$ for all $u\in(1,2)$ and $\chi(u)=1$ for all $u\geq2.$
Assume moreover that $\chi$ is non-decreasing in $\mathbb{R}.$ For
all $\theta\in(0,1/2),$ let $\chi_{\theta}$ be the function defined
by
$$\forall \,u\in\mathbb{R},~~\chi_{\theta}(u)=\chi(u/\theta).$$
The function $\chi_{\theta}$ is such that $0\leq\chi_{\theta}\leq1,$
$0<\chi_{\theta}<1$ in $(-\infty,\theta],$ $0<\chi_{\theta}<1$ in
$(\theta,2\theta)$ and $\chi_{\theta}=1$ in $[2\theta,+\infty).$
Furthermore, the functions $\chi_{\theta}$ are
non-increasing with respect to $\theta,$ namely,
$\chi_{\theta_1}\geq\chi_{\theta_2}$ if
$\displaystyle{0<\theta_1\leq\theta_2<1/2.}$

We set $$f_{\theta}(x,y,u)=f(x,y,u)\,\chi_{\theta}(u)~~\hbox{for
all}~(x,y,u)\in\overline{\Omega}\times\mathbb{R}.$$ In other words,
we cut off the source term $f$ near $u=0.$

For each $\theta\in(0,1/2),$ the function $f_\theta$ is a
nonlinearity of  ``combustion'' type that satisfies
(\ref{f1}-\ref{f2}) and (\ref{combustion}) with the ignition
temperature $\theta.$ Therefore, Theorem \ref{Berest Hamel
combustion Thm} yields that the existence of a classical solution
$(c_{\theta}, u_{\theta})$ of (\ref{front}) with the nonlinearity
$f_{\theta}.$ Furthermore, the function $u_{\theta}$ is increasing
in $t$ and unique up to translation in $t$ and the speed
$c_{\theta}$ is unique and positive.

It was proved, through Lemma 6.1 and Lemma 6.2 in Berestycki, Hamel
\cite{BH1}, that the speeds $c_\theta$ are non-increasing with
respect to $\theta$ and
$$c_\theta\nearrow c^{*}(e)~~\hbox{as}~\theta\searrow0.$$
Consider a sequence $\theta_{n}\searrow0.$ Then, there exists
$n_0\in\mathbb{N}$ such that $c_{\theta_n}\geq c^{*}(e)-\delta$ for
all $n\geq n_0$ (or equivalently
$\displaystyle{\theta_n\leq\theta_{n_{0}}}$).

 In what follows, we fix
$\theta$ such that $\displaystyle{\theta<\theta_{n_0}}.$ One
consequently gets $c_{\theta}\geq c^{*}(e)-\delta.$ On the other
hand, it follows, from the construction of $f_{\theta},$ that $f\geq
f_{\theta}$ in $\overline{\Omega}\times\mathbb{R}.$ Together with
(\ref{strict ineq}), one obtains
\begin{eqnarray}
c_{\theta}>&\displaystyle{\inf_{\varphi\in
E}\sup_{(s,x,y)\in\mathbb{R}\times\overline{\Omega}}\frac{F[\varphi](s,x,y)+q\cdot\nabla_{x,y}\varphi+f_{\theta}(x,y,\varphi)}{\partial_{s}\varphi(s,x,y)}+q\cdot\tilde{e}}.
\end{eqnarray}
Thus, there exists a function $\psi\in E$ such that
\begin{equation}\label{c_theta > strictly}
c_{\theta}>\displaystyle{\frac{F[\psi](s,x,y)+\,q\cdot\nabla_{x,y}\psi(s,x,y)+f_{\theta}(x,y,\psi)}{\partial_{s}\psi(s,x,y)}+q(x,y)\cdot\tilde{e}}.
\end{equation}
However, $\psi_{s}(s,x,y)>0$ for all
$(s,x,y)\in\mathbb{R}\times\overline{\Omega}.$ Thus, the inequality
(\ref{c_theta
> strictly}) can be rewritten as
\begin{equation}\label{Lpsi+f theta}
L\psi(s,x,y)+f_{\theta}(x,y,\psi)<0
~~\hbox{in}~\mathbb{R}\times\Omega,
\end{equation}
with $\psi\in E$ and $L$ is the operator defined in (\ref{degen
elliptic eq sat by phi}) for $c=c_\theta.$

For each $\tau\in\mathbb{R},$ we define the function $\psi^{\tau}$
by $$\psi^{\tau}(s,x,y)=\psi(s+\tau,x,y)~\hbox{for
all}~(s,x,y)\in\mathbb{R}\times\overline{\Omega}.$$ Since the
coefficients of $L$ are independent of $s,$ the later inequality also
holds for all functions $\psi^{\tau}$ with $\tau\in\mathbb{R}.$ That
is,
\begin{equation}\label{Lpsi tau + f theta <}
L\psi^{\tau}(s,x,y)+f_{\theta}(x,y,\psi^{\tau})<0
~~\hbox{in}~\mathbb{R}\times\overline{\Omega}.
\end{equation}

\vskip0.4cm
 \textit{Step 3.} For the fixed $\theta$ (in step 2), the function
$f_{\theta}$ is a ``combustion'' nonlinearity whose ignition
temperature is ${\theta}.$ There corresponds a solution
$(c_{\theta},u_{\theta})$ of (\ref{front}) within the nonlinear
source $f_{\theta}.$ We define $\phi_{\theta}$ by
$$\displaystyle{\phi_{\theta}(s,x,y)=u_{\theta}\left(\frac{s-x\cdot e}{c_{\theta}},x,y\right), ~\hbox{for all}~(s,x,y)\in\mathbb{R}\times\overline{\Omega}.}$$

Referring to section \ref{change of variables section}, one knows that $\phi_{\theta}\in E$
and thus it satisfies the following equation
\begin{equation}\label{Lphi theta =}
L\phi_\theta(s,x,y)+f_{\theta}(x,y,\phi_\theta)=0
~~\hbox{in}~\mathbb{R}\times\overline{\Omega}.
\end{equation}

Now, the situation is exactly the same as that in step 2 of the
proof of formula (\ref{min max combustion}) because the nonlinearity
$f_{\theta}$ is of ``combustion'' type. The little difference is
that $f$ (in step 2 of the proof of formula (\ref{min max
combustion})) is replaced here by $f_{\theta},$ and the function
$\phi$ of equation (\ref{Lphi+f(x,yphi)=0}) is replaced by the
function $\phi_{\theta}$ of (\ref{Lphi theta =}). Thus, following
the arguments of subsection \ref{proof of min max combustion} and using the same tools of ``step 2'' as in the
proof of formula (\ref{min max combustion}), one gets that the
(\ref{strict ineq}) is impossible and that completes the proof
of formula (\ref{min max ZFK }). \hfill $\Box$

\begin{remark} We found that one can use another argument (details are below) different from the sliding method in order to prove the $\min-\max$ formul{\ae} for the speeds of propagation whenever $f$ is a homogenous (i.e $f=f(u)$) nonlinearity of ``combustion'' or ``ZFK''type and $\Omega=\R^{N}.$ Meanwhile, the sliding method, that we used in the proofs of formul{\ae} (\ref{min max combustion}) and (\ref{min max ZFK }), is a unified argument that works in the general heterogenous periodic framework.
\vskip 0.2cm
\rm{\textbf{Another proof of formul{\ae} (\ref{min max combustion}) and (\ref{min max ZFK }) in a particular framework:}\\ Here, we assume that $f=f(u),$ and $\Omega=\R^{N}.$ Following the same procedure of ``step 1'' in the previous proof, one gets the inequality
$$c^{*}(e)\geq\displaystyle{\inf_{\varphi\in E}\sup_{(s,x,y)\in\mathbb{R}\times\overline{\Omega}}R\varphi(s,x,y)}.$$

Now, to prove the other sense of inequality, we assume that
$$c^{*}(e)>\displaystyle{\inf_{\varphi\in E}\sup_{(s,x,y)\in\mathbb{R}\times\overline{\Omega}}R\varphi(s,x,y)},$$
 and we assume that $f$ is of ``ZFK'' type \footnote{The case where $f$ is of ``combustion'' type follows in a similar way.}.
Then, as it was explained in ``step 2'' of the previous proof, one can find $\psi\in E,$ $\delta>0,$ $\theta>0,$ and $d>0$ such that $$c^{*}(e)-\delta<d<c_\theta<c^{*}(e)$$ where
$$
\forall(s,x,y)\in\mathbb{R}\times\overline{\Omega},~~d>c^{*}(e)-\delta>R\psi(s,x,y),$$
and
$f_{\theta}(u)=f(u)\,\chi_{\theta}(u) \leq f(u)\hbox{ for all } u\in\R$ is of ``combustion'' type ($c_\theta$ is the speed of propagation, in the direction of $-e,$ of pulsating travelling fronts solving (\ref{front}) with the nonlinearity $f_\theta$ and the domain $\Omega=\R^{N}$).

Hence, for all $(t,x,y)\in\R\times\R^{N},$
\begin{eqnarray}
d>&\displaystyle{\frac{F[\psi](s,x,y)+\,q\cdot\nabla_{x,y}\psi(s,x,y)+f_{\theta}(\psi)}{\partial_{s}\psi(s,x,y)}+q(x,y)\cdot\tilde{e}}.
\end{eqnarray}
Let $\tilde{u}(t,x,y)=\psi(x\cdot e+dt,x,y).$ As it was explained in section \ref{change of variables section}, the function
$\tilde{u}$ satisfies
\begin{eqnarray}\label{super sol}
    \left\{
      \begin{array}{ll}
\tilde{u}_t -\nabla\cdot(A(x,y)\nabla \tilde{u})-q(x,y)\cdot\nabla \tilde{u}-f_{\theta}(\tilde{u})>0,\;
t\,\in\,\mathbb{R},\;(x,y)\,\in\,\overline{\Omega},
\vspace{4 pt}\\
        \nu \cdot A\;\nabla \tilde{u}(t,x,y) =0,\;
        t\,\in\,\mathbb{R},\;(x,y)\,\in\,\partial\Omega,\vspace{4 pt}\\
 \displaystyle{\forall\, k\in\prod^{d}_{i=1}L_i\mathbb{Z},\; \forall\,(t,x,y)\in\mathbb{R}\times\overline{\Omega}},
\;\displaystyle{\tilde{u}(t+\frac{k\cdot e}{d},x,y)=\tilde{u}(t,x+k,y)}  \hbox{,} \vspace{4 pt}\\
          0\,\leq\,\tilde{u}\leq\,1.
      \end{array}
    \right.
\end{eqnarray}

Let $0\leq u_0(x,y)\leq 1$ be a function in $C(\R^{N})$ such that $u_0(x,y)\rightarrow0$ as $x\cdot e\rightarrow-\infty,$ and
$u_0(x,y)\rightarrow1$ as $x\cdot e\rightarrow +\infty,$ uniformly in $y$ and all directions of $\R^{d}$ which are orthogonal to $e.$ Let $u$ be a pulsating front propagating in the direction of $-e$ with the speed $c_\theta$ and solving the initial data problem
\begin{eqnarray}\label{spread}
    \left\{
      \begin{array}{ll}
u_t =\nabla\cdot(A(x,y)\nabla u)+q(x,y)\cdot\nabla u+\,f_\theta(u),\;
t>0,\;(x,y)\,\in\,\overline{\Omega},
 \vspace{3 pt}\\
u(0,x,y)=u_{0}(x,y),\vspace{3 pt}\\
        \nu \cdot A\;\nabla u(t,x,y) =0,\;
        t\,\in\,\mathbb{R},\;(x,y)\,\in\,\partial\Omega.
        \end{array}\right.
        \end{eqnarray}
Having $f_\theta(u)$ as a ``combustion'' nonlinearity, it follows from J. Xin \cite{JXin} (Theorem 3.5) and Weinberger \cite{weinberger}, that
\begin{equation}\label{spreading}
\begin{array}{c}
\forall r>0,~~\ds{\lim_{t\rightarrow+\infty}\ds{\sup_{|x|\leq r}u(t,x-cte,y)=0}}\hbox{ uniformly in $y,$ for every } c>c_\theta, \\
\hbox{and  }\ds{\lim_{t\rightarrow+\infty}\ds{\inf_{|x|\leq r}u(t,x-cte,y)=1}}\hbox{ uniformly in $y,$ for every } c< c_\theta.
\end{array}
\end{equation}

This means that the speed of propagation $c_\theta$ corresponding to (\ref{front}) is equal to the spreading speed in the direction of $-e$ when the nonlinearity is of ``combustion'' type and the initial data $u_0$ satisfies the above conditions.

For all $(t,x,y)\in[0,+\infty)\times\o,$ let $w(t,x,y)=\tilde{u}(t,x,y)-u(t,x,y).$ It follows, from (\ref{super sol}) and (\ref{spread}), that
  \begin{eqnarray}\label{parabolic on w}
    \left\{
      \begin{array}{ll}
w_t -\nabla\cdot(A(x,y)\nabla w)-q(x,y)\cdot\nabla w+\,bw>0,\;
t>0,\,(x,y)\,\in\,\overline{\Omega},\vspace{4 pt}\\
\forall(x,y)\in\o,~w(0,x,y)\geq0,\vspace{4 pt}\\
        \nu \cdot A\;\nabla w(t,x,y) =0,\;
        t\,\in\,\mathbb{R},\;(x,y)\,\in\,\partial\Omega,
        \end{array}\right.
        \end{eqnarray}
        for some $b=b(t,x,y)\in C(\ro).$ The parabolic maximum principle implies that $w\geq0$ in $[0,+\infty)\times \o.$ In other words,
$$\forall (t,x,y)\in [0,+\infty) \times\o,~~u(t,x,y)\leq \tilde{u}(t,x,y).$$
  However, for all $c>d,$ $$\ds{\lim_{t\rightarrow+\infty}\tilde{u}(t,x-cte,y)=\lim_{t\rightarrow+\infty}\psi(x\cdot e+(d-c)t,x-cte,y)=0}$$ locally in $x$ and uniformly in $y$ (since $\psi\in E$). Consequently,
$$\ds{\forall c> d,~\forall r>0,~~\lim_{t\rightarrow+\infty}\ds{\sup_{|x|\leq r}u(t,x-cte,y)=0}}\hbox{ uniformly in $y$}.$$
Referring to (\ref{spreading}), one concludes that $d\geq c_\theta$ which is impossible ($d<c_\theta$).
Therefore, our assumption that $c^{*}(e)>\displaystyle{\inf_{\varphi\in E}\sup_{(s,x,y)\in\mathbb{R}\times\overline{\Omega}}R\varphi(s,x,y)}$ is false
and that completes the proof of formula (\ref{min max ZFK }) in the case where $f=f(u)$ and $\Omega=\R^{N}.$\hfill$\Box$
}
\end{remark}

\section*{Acknowledgements}
I am deeply grateful to Professor Fran\c{c}ois Hamel for his valuable directions and advices during the preparation of this paper. I would like also to thank Professor Andrej Zlato\v{s} for his important comments, given while reading a preprint of this work, which allowed me to consider a wider family of ``ZFK'' nonlinearities in Theorem \ref{thm min max ZFK} and thus get a more general result.

\end{document}